\documentclass[a4paper,11pt]{amsart}

\usepackage[centertags]{amsmath}
\usepackage{amsfonts}
\usepackage{amssymb}
\usepackage{amsthm}
\usepackage{newlfont}
\usepackage{euscript}
\usepackage{wrapfig}
\usepackage{graphicx}

\theoremstyle{plain}
\newtheorem{thm}{Theorem}

\newtheorem{lem}{Lemma}
\newtheorem{prop}{Proposition}

\theoremstyle{definition}

\theoremstyle{remark}
\theoremstyle{remark}
\newtheorem{rem}{Remark}[subsection]

\numberwithin{equation}{subsection}

\renewcommand{\O}{{\EuScript O}}

\newcommand{\eps}{\varepsilon}

\newcommand{\s}{\sigma}

\newcommand{\jj}{\mathbf{i}}
\newcommand{\rr}{\mathbf{r}}

\newcommand{\1}{1\!\!{\mathrm I}}

\renewcommand{\P}{{\EuScript P}}

\newcommand{\G}{{\EuScript G}}

\newcommand{\Deu}{{\EuScript D}}

\renewcommand{\Re}{{\Bbb R}}

\newcommand{\ax}{\Re^+}

\newcommand{\toP}{\mathop{\longrightarrow}\limits^{P}}

\newcommand{\cth}{\mathop{\hbox{cth}}\nolimits}
\newcommand{\sh}{\mathop{\hbox{sh}}\nolimits}
\newcommand{\ch}{\mathop{\hbox{ch}}\nolimits}
\newcommand{\arcth}{\mathop{\hbox{arcth}}\nolimits}
\renewcommand{\th}{\mathop{\hbox{th}}\nolimits}

\begin{document}

\title{A limit theorem for diffusions on graphs with variable configuration}
\author{Alexey M. Kulik}%

\abstract{A limit theorem for a sequence of diffusion processes on
graphs is proved in a case when vary both parameters of the
processes (the drift and diffusion coefficients on every edge and
the asymmetry coefficients in every vertex), and configuration of
graphs, where the processes are set on. The explicit formulae for
the parameters of asymmetry for the vertices of the limiting graph
are given in the case, when, in the pre-limiting graphs, some
groups of vertices form knots contracting into a points.}
\endabstract

\date{}
\subjclass[2000]{60J60, 60F17, 60H10}%
\keywords{Diffusion process on a graph, Walsh's Brownian motion,
graphs with variable configuration}
\address{Institute of Mathematics,
Ukrai\-ni\-an National Academy of Sciences, 3, Tereshchenkivska
Str., Kyiv 01601, Ukraine}
 \maketitle

\centerline{\textsc{Introduction}} \vskip 20pt The main object
considered in the present paper is diffusion processes on graphs;
a graph is treated as a one-dimensional topological space with
branching points, rather than as a discrete scheme. Such processes
arise naturally, on the one hand, in a number of applied models
(e.g., in a model describing the motion of nutrients  in the root
system of a plant, see \cite{FrDU}) and, on the other hand, in
some theoretical constructions (e.g., in the study of small random
perturbations of Hamiltonian dynamical systems, see \cite{Fre}, or
in the study of the asymptotic geometric properties of discrete
groups, see \cite{V}). Such processes possess a number of
interesting and nontrivial internal structural peculiarities. Let
us mention one of them, that was revealed by B.S.Tsirelson and
consists in the following (see detailed review in
\cite{Tsirelson}): the typical diffusion on a graph generates a
filtration, that cannot be obtained from the filtration generated
by some (possibly, infinite-dimensional) Brownian motion, "in real
time", i.e., by means of a {\em morphism}.

The full description of a diffusion process on a graph in terms of
its infinitesimal operator is given in \cite{FreWen}. In a
nonformal way, such a process can be described as a mixture of the
motions "along an edge"\phantom{} and "in a neighborhood of a
vertex". A motion of the first type is described by a
one-dimensional diffusion process and is defined by its
coefficients of drift and diffusion. To describe a motion of the
second type, it is necessary to set additionally  the parameters
playing the role of boundary conditions at a vertex, that define
the behavior of a process in the vertex. In \cite{FreWen}, these
objects are called "the gluing parameters". We call them also "the
asymmetry parameters", since the construction of the process, that
we use as the  basic one, differs from the one developed in
\cite{FreWen}. The purpose of this paper consists in the
description of the limit behavior of a diffusion process in the
situation where both the above-mentioned parameters of the process
and the graph itself vary.  If the limiting lengths of edges of
the graph are nonzero (i.e., the vertices of the graph do not
"glue together", and the configuration of the graph, in fact, does
not vary), then the required limiting result is very similar to
the standard limit theorems of the theory of diffusion processes.
The new specific problems requiring a separate analysis arise in
the case where some groups of vertices form knots contracting into
a points. In this case, the description of the limiting behavior
of diffusions can be investigated in the framework of the approach
developed in \cite{FreWen}, and non-trivial problem consists in
calculation of the asymmetry parameters of the limiting process.
As an example of the expressions that can be finally obtained, we
give the model, in which the graph consists of two vertices and
three edges: the edge with length $\eps$ joining the vertices and
two half-lines beginning at these vertices. If the diffusion on
each edge is a Brownian motion, then its distribution converges,
as $\eps\to 0$, to the distribution of a skew Brown motion on a
real line (\cite{Zaitseva}). Its asymmetry parameters are set by a
single "skewing parameter" $q$ ($p_{\pm}={1\pm q\over 2}$) which
is calculated in the given model as follows: if $q_1,q_2$ are
analogous parameters for the vertices of a prelimiting graph, then
$q=\th (\lambda_1+\lambda_2)$, where $\lambda_{1,2}=\arcth
q_{1,2}$, $\th c\equiv {e^c-e^{-c}\over e^c+e^{-c}}$ is a
hyperbolic tangent, and $\arcth$ is the function inverse to $\th$.
This simple example shows that the determination of the asymmetry
parameters for a vertex, that is obtained as a result of the
contraction of a group of vertices into a single point, is a
nontrivial problem. In the present paper, we propose a method of
the solution of this problem.

Let us describe one possible application of the main result of the
paper. Earlier, we have mentioned the Tsirelson's result that
states that if, for some vertex of a graph, at least 3 it's
asymmetry parameters are nonzero (i.e.,  the vertex of a graph has
multiplicity $\geq 3$ and is \emph{a triple point}), then the
diffusion $X$ is essentially singular in the following sense:
there does not exist any morphism of the filtration, generated by
any Brownian motion, to the filtration  generated by $X$
(\cite{Tsirelson}). We call further the singularity, related to
the presence of a triple point (in the above-mentioned sense), the
\emph{Tsirelson's singularity}.

The theory of diffusion  processes includes a number of results,
in which the diffusions, containing singularities of some type
(such as an asymmetric semipermeable boundary at some point, or a
partial reflection with delay at the boundary of a domain), are
represented as the weak limits of nonsingular diffusions (see
\cite{Kulinich},\cite{Chitashvili}). In these results, typically,
the coefficients of prelimiting diffusions, in a certain sense,
model the singular terms that are present in the limiting process
(such as the singular drift coefficient $a=q\delta_0$ of a skew
Brownian motion). The Tsirelson's singularity is related to the
structure of the phase space (the presence of a nontrivial
branching point), rather than to properties of the coefficients of
the process. The main theorem of the present paper allows us to
represent a process, possessing the Tsirelson's singularity, as
the weak limit of the processes without such singularities. We
will construct such a representation, changing the phase space of
the process and representing a vertex with multiplicity $\geq 3$
as a result of "the contraction into a point"\phantom{} of a knot,
whose vertices have multiplicity $\leq 2$ (see Example 2 below).
Such a trick gives, in perspective, the possibility to study the
properties of stochastic flows (i.e., the processes describing a
motion of a families of points) corresponding to the diffusion on
a graph, that possesses the Tsirelson's singularity and for which,
by this reason, one cannot define the common law of motion of the
family of points as the strong solution to a system of SDE's.

\vskip 20pt \centerline{\textsc{1. Basic notation and
constructions}} \vskip 20pt

{\bf 1.1. Phase space.} Everywhere in what follows, a graph means
a connected metric space composed of a finite number of subsets
(edges of the graph) homeomorphic to a segment or a half-line. It
is assumed that the only intersection points for these subsets are
the images of the ends of a segment or a half-line; all such
points are the vertices of the graph.

It is convenient to consider the graph to be oriented, by assuming
that if a point, whose motion is described by a diffusion process,
is not positioned on a definite edge, then this point can move
into this edge, only by passing through the separated vertex of
the edge (its "beginning"). It is clear that this does not
restrict the generality, because any non-oriented
edge-"segment"\phantom{} can be represented as two oriented
copies; for an edge-"ray"\phantom{}, its beginning, obviously, is
the single vertex belonging to it. For the vertex $i$ of a graph
$\G$, we denote, by $\P_i$, the family of vertices joined with it
by edges beginning at the vertex $i$, and these edges are denoted
as $L_{i,j}^r,j\in \P_i,r=1,\dots, R(i,j)$. The necessity to
introduce the additional parameter $r$ is caused by that two
vertices can be joined by several edges. Analogously, the edges
that leave the vertex $i$ and are homeomorphic to a half-line are
denoted as $L_{i,\infty}^r, r=1,\dots, R(i,\infty)$. The interior
of the set $\left[\bigcup_{j\in\P_i}\bigcup_r
L_{i,j}^r\right]\cup\left[\bigcup_r L_{i,\infty}^r\right]$ is
denoted as $\O_i$ and called the maximum neighborhood of the
vertex $i$.

On each edge, we introduce a natural parametrization in the
following way: for the edge $L_{i,\infty}^r$, we consider the
homeomorphism with $[0,+\infty)$ to be fixed and define the
coordinate of a point on the edge as the number corresponding to
it via this homeomorphism. The coordinates of points on the edges
$L_{i,j}^r,j\in \P_i$ are determined analogously, by considering
the homeomorphism of the edge with $[0,l_{i,j}^r]$ to be fixed.

We make no assumptions that the graph is imbedded in any larger
metric space (e.g., that it is planar). On the other hand, the
natural parametrization allows us to homeomorphically imbed each
of the maximum neighborhoods $\O_i$ in $\Re^2$, the corresponding
image being a part of a bundle of half-lines in $\Re^2$. If this
does not cause misunderstanding, we will omit the corresponding
homeomorphism, by considering that any $\O_i$ is a part of a
bundle of half-lines in $\Re^2$ and that any edge is a part of
$\Re^+$.

{\bf 1.2. Construction of the process.} In our consideration, we
will use two constructions of a diffusion process on a graph. One
construction sets the process in terms of its infinitesimal
operator and is given in \cite{FreWen}. The other, more explicit
construction describes the diffusion process in terms of its
excursions. We cannot give a reference, where the required
construction would be described in  the generality sufficient for
our purposes, therefore, we give its description here. Of course,
we do not pretend for a priority, because various versions of such
a construction were given earlier (see \cite{ItoMcKean}, \S4.2 and
\cite{W},\cite{BPY}), and its main idea is widely known. In the
next subsection, we will show that both versions of the
construction of a diffusion process on a graph lead to the same
object.

First, let a graph $\G$ be a bundle of half-lines, i.e., it has
one vertex and $R$ edges-"rays". The \emph{Walsh's Brownian
motion} is set on such a graph as follows (\cite{W},\cite{BPY}):
it is a continuous Feller Markov process; on the start from a
point on one of the rays, up to the first moment to hit of the
vertex, this process propagates along this ray as the Wiener
process. On the start from the vertex, the process is constructed
as follows: one take the Brownian motion $B$ on $[0,+\infty)$ with
reflection at zero, and construct a sequence of identically
distributed random values $\{\eps_n,n\geq 1\},$ that  are
independent of the process $B$ and one another and such that
$P(\eps_1=k)=p_k, k=1,\dots,R$, where $p_1,\dots,p_R$ are the
given numbers ("the asymmetry parameters"\phantom{} of the
vertex), $p_1+\dots+p_R=1$. Then one enumerate, in an arbitrary
measurable way, all the excursions $\{[a_m,b_m]\}$ of the process
$B$ at zero and denote,  by $m(t)$ for an arbitrary time moment
$t$, the random variable setting the number of a current
excursion, i.e., such natural number that $a_{m(t)}<t<b_{m(t)}$.
Such a variable is well defined on the set $\{B_t>0\}$ having
probability 1. Now the value of the required process $X$ at a time
moment $t$ is determined in the following way: it is located on a
ray with the number $\eps_{m(t)}$, and its coordinate on this ray
is equal to $B_t$.

Next, let the graph $\G$ be the same one as above, and let, for
each edge $L^r$, the functions $a^r,\sigma^r$ on the ray
$(0,+\infty)$, which are interpreted as the coefficients of drift
and diffusion of the process on the edge, be set. We assume that
these functions are measurable, the functions $\sigma^r,
[\sigma^r]^{-1}$ are locally bounded, and $a^r$ is locally
integrable. In addition, in order to shorten the consideration and
to exclude the possibility of an "explosion", we assume that the
functions $a^r,\sigma^r$ possess at most linear growth at
infinity. We introduce a new parametrization
$$\hat x= S^r(x)\equiv \int_0^x
\exp[-\int_0^y{2a^r(z)\over [\sigma^r(z)]^2}\,dz]\,dy,\quad x\in
L^r,\eqno(1) $$ and, for the given collection
$\wp=\{p_1,\dots,p_R\}$ of the asymmetry parameters of a vertex,
construct a Walsh's Brownian motion $\hat X$ with these
parameters. We introduce a process $\{r(s),s\geq 0\}$ such that
$\hat X_s\in L^{r(s)}, s\geq 0,$ almost surely and put
$$\theta_t=\int_0^t
\left\{\sigma^{r(s)}(\hat X_s)\cdot [S^{r(s)}]'(\hat X_s)
\right\}^2ds,\quad \tau_t=\theta_t^{-1}\equiv \inf\{u|\theta_u\geq
t\}, \quad t\geq 0.\eqno(2)
$$
Now, we put $X_t=S^{-1}(\hat X_{\tau_t}) ,t\geq 0$, where $S^{-1}$
is the change of a parametrization on the initial graph which is
inverse to $S$. By construction, $X$ is the Feller Markov process
with continuous trajectories. The proof of this assertion is
exactly the same as that for an analogous proposition about the
general one-dimensional diffusion (see \cite{ItoMcKean}, Chap. 3).
We interpret process $X$ as the diffusion process on a bundle of
half-lines with the coefficients of drift and diffusion on rays
$\{a^r,\sigma^r\}$ and the asymmetry parameters of the vertex
$\wp$.

At last, let a graph $\G$ be arbitrary, let the collection of
asymmetry parameters $\wp_i$ be set for its each vertex, and let
the coefficients of drift and diffusion
$\{a_{i,j}^r,\sigma_{i,j}^r\}$ be set for each edge $L_{i,j}^r$.
Supposing, e.g., that
$a_{i,j}^r(x)=\sigma_{i,j}^r(x)=1,j\not=\infty,x\geq l_{i,j}^r$
(or $a_{i,j}^r(x)=a_{i,j}^r(l_{i,j}^r),
\sigma_{i,j}^r(x)=\s_{i,j}^r(l_{i,j}^r),j\not=\infty,x\geq
l_{i,j}^r$, if $a_{i,j}^r,\s_{i,j}^r$ are continuous on
$[0,l_{i,j}^r]$), we can consider that these functions are given
on the whole $\Re^+$. Now,  the Feller Markov process with
continuous trajectories is well defined on $\G$ via the following
convention: on the start from a point lying in the maximum
neighborhood $\O_i$ of any vertex $i$, it moves as a diffusion
process on the bundle of rays with the parameters $\wp_i,\{
a_{i,j}^r,\sigma_{i,j}^r,j\in\{1,\dots,N\}\cup\{\infty\},r=1,\dots,R(i,j)\}$
till the exit  from $\O_i$. At the moment of exit from $\O_i$, it
is located at some vertex $j\in \P_i$, and, after this moment and
up to the next exit time (from $\O_j$), it moves as a diffusion
process on the bundle of rays in the maximum neighborhood $\O_j,$
and so on. The process, constructed in the way described above, is
the diffusion process on $\G$ with the parameters $\{\wp_i,
a_{i,j}^r,\sigma_{i,j}^r\}$.

Let us note that the constructed process spends zero time at every
point of the phase space with probability 1. A wider class of the
processes with "sticky"\phantom{} points can be constructed by
means of a time change (see the details, e.g., in \cite{portenko},
Section 3.3). In order to  shorten the exposition,  we exclude
such processes from consideration.

{\bf 1.3. Martingale description and infinitesimal
characteristics of the process.}

To describe the infinitesimal characteristics of the process
constructed in the previous section, it is sufficient to consider
the case where the graph $\G$ is a bundle of half-lines. For such
$\G$, let the diffusion $X$ with the parameters
$\wp=\{p_r\},\{a^r,\sigma^r\}$ be given. We assume that the
functions $a^r,\sigma^r$, in addition to the above-imposed
conditions, are continuous on $[0,+\infty)$. By $X^r$, we denote a
coordinate process on the $r$-th edge: $X^r_t$ is equal to the
coordinate $X_t$ on the edge $L^r$, if $X_t$ lies on this edge at
a time moment $t$, and $X^r_t$ equals zero otherwise. A vertex of
the graph is denoted by the letter $O$. The process $X$ can be
described in terms analogous to those of the Skorokhod problem for
a process on a half-line with reflection at a point of the
boundary.

\begin{prop}\label{mart_problem_prop} There exists a nondecreasing process
$V_t$ such that

(i) it increases only  when the vertex is visited by the process
$X$, i.e.,
$$
\int_0^t\1_{\{X_s\not=O\}}\,dV_s=0\quad \hbox{a.s.,}\quad t\geq 0;
$$

(ii) for an arbitrary edge $L^r$, the process
$$
M^r_t\equiv X^r_t-p_rV_t-\int_0^ta^r(X_s^r)\1_{\{X_s\in
L^r\}}\,ds$$ is a continuous martingale with the quadratic
variation
$$
\langle M^r\rangle_t=\int_0^t[\sigma^r(X_s^r)]^2\1_{\{X_s\in
L^r\}}\,ds.$$
\end{prop}

\emph{Proof.} For a Walsh's Brownian motion, this assertion
follows from the reasoning analogous to that given in
\cite{Tsirelson}, Section 3: each of the processes $X^r$ belongs
to the class $\Sigma_+$ (\cite{RY}, VI.4.4.), i.e., can be
represented as a sum $M^r+V^r$, where $M^r$ is a local martingale,
and $V^r$ is a nondecreasing process such that
$\int_0^t\1_{\{X^r_s>0\}}dV^r_s\equiv 0$. It is easy to verify
that the quadratic characteristic of $M^r$ is equal to
$\int_0^t\1_{\{X_s^r>0\}}\,ds$. On the other hand (see \cite{RY},
VI.4.4.),
$$
V^r_t=\lim_{\eps\to 0+}{1\over
2\eps}\int_0^t\1_{\{X_s^r\in(0,\eps)\}}\,ds.
$$
Let $V_t$ be the local time of the process $X$ at a vertex $O$, i.e.,
$$V_t=\lim_{\eps\to 0+}{1\over
2\eps} \int_0^t\1_{\{dist(X_s,O)<\eps\}}\,ds.
$$
It is easy to see that both $V^r$ and $V$ are  $W$-functionals of
the process $X$. Moreover, their characteristics $f^r$ and $f$ are
connected by the relation $f^r=p^r\cdot f$ by construction. Using
Theorem 6.3 in \cite{dynkin}, we get $V^r=p^r\cdot V$. That is,
Proposition 1 is valid for a Walsh's Brownian motion, and the
process $V$ is the local time of the process $X$ at zero.

Now let $X_t=S^{-1}(\hat X_{\tau_t})$, where  $\hat X$ is  a
Walsh's Brownian motion (we use the construction of the previous
section), and $\hat V_t$ is its local time at the vertex. Then
Proposition 1 is valid with $V_t=\hat V_{\tau_t}$, that  follows
from the It\^o formula (note that $[S^r]'(0)=1$).

Assertions (i),(ii) can be naturally interpreted as a version of
the martingale problem for the pair of processes $(X,V)$. Another
version of the martingale problem for the process $X$ is given
below.

\begin{prop}\label{mart_problem_prop2} Let a continuous function
$\phi$ on $\G$ be equal, on every edge $L^r$, to some  function
$\phi^r\in C^2([0,+\infty)$. Denote
$A\phi(x)=a^r(x)[\phi^r]'(x)+{1\over 2}[\sigma^r(x)]^2
[\phi^r]''(x), x\in L^r$ and
$\Delta_O(\phi)=\sum_rp_r[\phi^r]'(0)$.

Then, for an arbitrary $\phi$ satisfying the above-indicated
condition and such that $\Delta_O(\phi)=0$, the process
$$
M_t^\phi\equiv \phi(X_t)-\int_0^tA\phi(X_s)\,ds
$$
is a continuous martingale.
\end{prop}

This assertion follows immediately via the It\^o formula from
Proposition 1 and the fact that the process $X$ spends  zero time
at $O$.

It is important that, as the following theorem shows, both the
martingale problems given in Propositions 1,2 are well-posed. By
$\Deu_A$, denote the set of continuous bounded functions $\phi$ on
$\G$, such that the  function $A\phi$ is well defined, continuous,
and bounded on $\G$, and the condition $\Delta_O(\phi)=0$ holds.

\begin{thm}\label{mart_problem_thm} 1) The operator $(A,\Deu_A)$
is an infinitesimal operator of the process $X$ constructed in the
previous subsection.

2) The process $X$ is the unique solution of the martingale
problem posed in Proposition 2, endowed by the given initial
distribution $P(X_0\in\cdot)$.

3) Let $V$ be the process constructed in the proof of Proposition
1. Then the pair of processes $(X,V)$ is the unique solution of
the martingale problem posed in Proposition 1, endowed by the
given initial distribution $P(X_0\in\cdot)$, that satisfies the
following conditions:

(i) $P(V_0=0)=1$;

(ii) with probability 1, the process $X$ spends zero time at the
vertex $O$.
\end{thm}

{\it Proof.} The fact that $(A,\Deu_A)$ is a pre-generator of the
process $X$ follows from Proposition 2. Theorem 3.1 in
\cite{FreWen} ensures the fact that $\Deu_A$ is the whole domain
of definition of the generator of the process $X$. The second
assertion follows from Theorems 3.1 and 2.2 in \cite{FreWen}. The
third assertion is a consequence of the second one and the It\^o
formula.

According to Theorem 1, the constructive description of a
diffusion process on a graph, presented in the previous section,
and the semigroup description presented in Section 3 in
\cite{FreWen} are equivalent.

{\bf 1.4. Limit theorem for a graph with constant configuration.}
The construction of Section 1.2 yields directly the following
limiting result for the sequence of diffusion processes $\{X^n\}$
on a graph $\G$ with the parameters
$\{\wp_i^n,a_{i,j}^{r,n},\s_{i,j}^{r,n}\}$.

\begin{thm}\label{lim_thm simple} Let
$$
a_{i,j}^{r,n}\to a_{i,j}^{r}, \hbox{ в } L_{1,loc}(\Re^+), \quad
\s_{i,j}^{r,n} \to \s_{i,j}^{r}
$$
locally uniformly on $[0,+\infty)$, and let
$$
\wp_i^n=\{p_{i,j}^{r,n},j\in\P_i\cup\{\infty\}, r=1,\dots
R(i,j)\}\to \wp_i=\{p_{i,j}^{r},j\in\P_i\cup\{\infty\}, r=1,\dots
R(i,j)\}
$$
componentwise (recall that we suppose that the coefficients
$a_{i,j}^{r,n},\s_{i,j}^{r,n}$ are given on $[0,+\infty),$ by
setting them, if necessary, to a constant on
$(l_{i,j}^{r,n},+\infty)$).

Then the sequence of the distributions of processes $\{X^n\}$ converges weakly in
$C(\Re^+,\G)$ to the distribution of a diffusion
process $X$ with the parameters $\{\wp_i,a_{i,j}^{r},\s_{i,j}^{r}\}$.
\end{thm}

{\it Proof.} In each of the maximum neighborhoods $\O_i$, we
consider the changes of the phase variable $S^n,[S^n]^{-1},$ that
are defined by the coefficients $X^n$ on the corresponding edges.
By virtue of the imposed conditions, $S^n,[S^n]^{-1}$ converge
locally uniformly to $S,[S]^{-1}$,  defined by the coefficients of
$X$. By performing the change of a phase variable for $X^n$ in
each of the neighborhoods $\O_i$, that  is inverse to (1),  and
the change of the time that is inverse to (2), we obtain the
process $\hat X^n$ being a composition of the Walsh's Brownian
motions, switching themselves at the time moment of the transition
from one neighborhood to another one. The integrands $\{\sigma^n
[S^n]'\}^2$ in the time change (2) time also converge on each edge
locally uniformly to $\{\sigma [S]'\}^2$. The probability of the
event, that at least $M$ transitions from one neighborhood to
another one occur for the process $\hat X^n$, can be estimated
uniformly in both $n$ and the starting point by a term  of the
form $C\alpha^M$,  where $\alpha\in(0,1),C$ are some constants
(this follows from the estimate $P(\sup_{t\leq \eps}|W(t)|>c)\leq
e^{-{Kc^2\over \eps}}$ for the Wiener process $W$; see also an
analogous estimate for the It\^o processes in \cite{WatIkeda},
Lemma 8.5). Thus, by virtue of the strong Markov property for
$X^n, X$, the proof of Theorem 2 is reduced to the following.
There exists the sequence of the Walsh's Brownian motions
$\{Z^n\}$ on a given bundle of rays starting from the vertex, and
the corresponding collections $\wp^n$ converge componentwise to
the collection $\wp$ for the Walsh's Brownian motion $Z$. There is
also a family of neighborhoods of the vertex $\O^n$, $\lim_n
\O^n=\O$, and $\tau_{\O^n}^{Z^n}$ and $\tau_{\O}^{Z}$ are the
moments of the exit of the processes from $\O^n$ and $\O$,
respectively. We have to show that the distributions of the pair
$(Z^n,\tau_{\O^n}^{Z^n})$ converge weakly to the distribution
$(Z,\tau_{\O}^{Z})$. This can be proved, by leaning on three
following simple assertions. First, $Z^n$ converges weakly to $Z$
(here, we can explicitly write the transition probabilities for
$Z^n, Z$). Secondly, $\tau_{\O}^{Z}$, for an arbitrary
neighborhood $\O$, is an almost surely continuous functional of
the trajectory $Z$. At last, $\tau_{\O}^{Z}$ is a monotonous and,
for almost all trajectories $Z$, continuous function of $\O$. Two
last facts are a consequence of analogous assertions for the
Brownian motion. Theorem 2 is proved.

 {\it Remarks. 1.} The assertion of Theorem 2 could be proved,
by using the general limit theorem 4.1 in  \cite{FreWen}. However,
the reasoning presented in the proof is important for us by
themselves, because it is, in essence, a part of the proof of the
main result of the present paper, Theorem 3.

{\it 2.} The assertion of Theorem 2 remains valid if we assume
that, for each $\{X^n\}$, its diffusion coefficients are written
w.r.t. its own parametrization $\Psi^n$, and $\Psi^n\circ
\Psi^{-1}\to \mathrm{id},  \Psi\circ [\Psi^n]^{-1}\to \mathrm{id}$
locally in $C^2(\Re^+)$ on each edge. Informally, this means that
the lengths of edges can vary, not tending to zero.

\vskip 20pt \centerline{\textsc{2. Limit theorem for graphs with
variable configuration}} \vskip 20pt

{\bf 2.1. Statement.} Further, we assume that the diffusion
processes $X^n$ are set on the graphs $\G^n$ with identical
combinatorial configuration (i.e., with the identical procedure to
join edges), but with different metrics. Formally, this
corresponds to the setting of different parameterizations $\Psi^n$
on the same graph $\G$. Nonformally, this means that the family of
edges is not changed, but their lengths vary.  We assume that the
characteristics $\{\wp_i^n,a_{i,j}^{r,n},\s_{i,j}^{r,n}\}$ of the
processes $X^n$ (relative to the corresponding parametrizations
$\Psi^n$) satisfy the conditions of Theorem 2 and consider the
question about the limiting behavior of the distributions of the
processes $X^n$.

We assume that the lengths $l_{i,j}^{r,n}$ of some edges
tend to zero. For $i,j$ such that $\exists r:
l_{i,j}^{r,n}\to 0,p_{i,j}^{r,n}\not \to 0$, we write $i
\dashrightarrow j$. For $i,j$ such that $\exists j_1,\dots, j_m:
i\dashrightarrow j_1, j_1 \dashrightarrow j_2,\dots,j_m
\dashrightarrow j$, we write $i \rightsquigarrow j$. We require that the
 following symmetry condition be satisfied:
$$
i \rightsquigarrow j \Leftrightarrow j \rightsquigarrow i.
\eqno(3)
$$
If condition (3) fails then  the  limiting process can fail to be
a diffusion.  Namely, the limiting process (in the sense of the
convergence of finite-dimensional distributions) can possess
discontinuous trajectories if condition (3) does not hold.

By construction, "$\rightsquigarrow$"\phantom{} is the equivalence
relation; the collection of vertices connected by this relation
will be called "a knot". It is natural to define a new graph $\hat
\G$, on which a limiting process will be defined at last, as a
graph, in which the vertices are the knots of the initial graph,
and the edges are those edges of the initial graph that have not
contracted into a point. We impose the following natural
condition:

\textbf{A.} On each of the edges $L_{i,j}^{r,n}$ not contracting
into a point, the parameterizations $\Psi^n$ converge in the sense
of the convergence of the $C^2$-diffeomorphisms of $\Re^+$ to a
certain parametrization $\Psi$, and $l_{i,j}^{r,n}\to
l_{i,j}^{r}>0$ as $n\to \infty$.

The phase space of the limiting process is the graph $\hat \G$
with the parametrization $\Psi$. To formulate the limiting result,
it is necessary to set the projection $\{\hat X^n\}$ of the
initial sequence $\{X^n\}$ on this space (see \cite{FreWen},
Section 4). This can be performed in the following way: if $X^n_t$
lies on a non-contracting edge of the graph $\G$, then $\hat
X^n_t$ lies on the corresponding edge of the graph $\hat \G$, and
its coordinate on this edge is obtained from the coordinate of
$X^n_t$ by the transformation $\Psi\circ [\Psi^n]^{-1}$. If
$X^n_t$ lies on a contracting edge, whose ends belong to the knot
$\hat i$, then $\hat X^n_t=\hat i.$ It is clear that the
trajectories of the process $\hat X^n$ are continuous.

The asymmetry parameters of the limiting process will be
determined by the internal structure of prelimiting knots, let us
introduce the necessary objects and assumptions. The main
assumption consists in that each knot is homogeneous, i.e., the
lengths of all internal edges in it tend to zero with the same
rate,
$$
\forall\,\hat i\,\,\,\exists\, \phi_{\hat i}(\cdot): \phi_{\hat
i}(n)\to 0\hbox{ and }\forall i,j\in \hat i, r=1,\dots,R(i,j)\quad
{l_{i,j}^{r,n}\over \phi_{\hat i}(n)}\to l_{i,j}^{r}>0,\, n\to
+\infty, \eqno(4)
$$
where the numbers $\{l_{i,j}^{r}\}$ are given. This assumption can
be weakened, but it cannot be removed at all. To shorten the
consideration, we omit the details here.  We note only that if
condition (4) is not imposed, we can faced with situations where
the limiting process spends a positive time at a knot with a
nonzero probability.

We recall that we assume that, for
each vertex $i$, the collections of the parameters $\wp_i^n$ converge to
a certain collection $\wp_i=\{p_{i,j}^{r},j\in\P_i\cup\{\infty\},
r=1,\dots R(i,j)\}$.

Let the knot $\hat i$ be fixed. We set $N_{\hat i}=\#\{i\in \hat
i\}$,
$$
\alpha^{\hat i}_{i,j}=\sum\limits_{r\leq R(i,j)} {p^{r}_{i,j}\over
l^{r}_{i,j}},\quad i,j\in\hat i,\quad \beta^{\hat i}_i=\sum_{j\in
\P_i\cap \hat i}\alpha^{\hat i}_{i,j},\quad i\in \hat i.
$$
Consider the $N_{\hat i}\times N_{\hat i}$-matrix $A^{\hat i}$
defined by
$$
A^{\hat i}_{i,j}={\alpha^{\hat i}_{i,j}\over \beta^{\hat
i}_i},\quad \quad i,j\in\hat i.
$$
The matrix $A^{\hat i}$ is the matrix of transition probabilities
for some Markov chain. By virtue of condition (3), all the states
of this chain form unique class of essential states. Therefore,
there exists the unique invariant distribution for the chain. We
denote this distribution by $\pi^{\hat i}$ and set the collection
$\hat \wp_{\hat i}$ in the following way: each edge $L_{\hat
i,\hat j}^{\hat r}$ of the graph $\hat \G$ is represented by some
edge $L_{i,j}^r$ of the graph $\G$ with $i\in \hat i,j\not \in
\hat i$ (it is possible that $\hat j=\infty$, then $j=\infty$).
For this edge,  we put
$$
\hat p_{\hat i,\hat j}^{\hat r}= P^{\hat i}\cdot {\pi^{\hat
i}_{i}\over \beta _i^{\hat i}}\cdot p_{i,j}^{r}, \eqno (5)
$$
where $P^{\hat i}$ is the normalizing factor which is defined by the
condition
$$
\sum_{\hat j\in\P_{\hat i}\cup \{\infty\},\hat r\leq R(\hat i,\hat
j)} \hat p_{\hat i,\hat j}^{\hat r}=1.
$$
At last, we set the functions $\hat a_{\hat i,\hat j}^{\hat
r},\hat \s_{\hat i,\hat j}^{\hat r}$ as the limits of the
functions $ a_{i,j}^{r,n},  \s_{i,j}^{r,n}$ on each edge $L_{\hat
i,\hat j}^{\hat r}$ corresponding to the edge $L_{i,j}^r$.  Now we
can formulate the main result of the present paper.

\begin{thm}\label{lim_thm strong} Let the
characteristics
$\{\wp_i^n,a_{i,j}^{r,n},\s_{i,j}^{r,n}\}$ of the processes $X^n$
satisfy the conditions of Theorem 2, and let
conditions (3) and (4) be satisfied. Let also the sequence of distributions
$\mu^n(\cdot)=P(\hat X^n(0)\in \cdot)$ converge weakly to some
measure $\mu$.
Then the sequence of distributions of the processes $\{\hat X^n\}$ in
$C(\Re^+,\hat \G)$ converges weakly to the distribution of the diffusion
process $\hat X$ with the above-set parameters $\{\hat \wp_{\hat i},\hat a_{\hat
i,\hat j}^{\hat r},\hat \s_{\hat i,\hat j}^{\hat r}\}$
and with the initial distribution $\mu$.
\end{thm}

{\bf 2.2. Proof.} The reasoning analogous to those used in the
proof of Theorem 2 allows us to restrict our  consideration to the
case where $a_{i,j}^{r,n}\equiv 0, \s_{i,j}^{r,n}\equiv 1$, all
the vertices of the initial graph $\G$ form a single knot, and
only nontrivial (i.e., not contracting into a point) edges are the
edges-"rays" that are homeomorphic to half-line.

\textbf{Remark 3.} For such a reduction, it is significant that
the conditions of convergence of the coefficients
$a^{r,n}_{i,j},\s^{r,n}_{i,j}$, which were formulated in Theorem
2, hold also for edges contracting into a point. Otherwise, the
assertion of Theorem 3 can be violated. For example, if for some
(not all) edges, contracting into a point, $a^{r,n}_{i,j}$ are
constant functions tending to $+\infty$ as $n\to +\infty$, then
the homogeneity condition (4) will be broken after a change of the
phase variable. If for some (not all) edges, contracting into a
point, $\s^{r,n}_{i,j}$ are constant functions, which tend
sufficiently rapidly to $0$ as $n\to +\infty$,  then the time
spent by the process $\hat X^n$ at the vertex $\hat i$, will not
tend to zero.

Let us proceed with the proof of Theorem 3 in the above-indicated
case. First of all, we note that the sequence of distributions of
the processes $\{\hat X^n\}$ is weakly compact in $C(\Re^+,\hat
\G)$. The simplest way to prove this,  is to use the criterion for
weak compactness given in \cite{Billingsley}, Theorem 8.2. For an
arbitrary $\eps>0$, the process $\{\hat X^n\}$ outside the vertex
neighborhood $B(O,\eps)$ with radius $\eps$ is the Brownian
motion. This  easily implies that, on every finite time interval
$[0,T]$, for the continuity modulus $w_T(\hat X^n,\delta) \equiv
\sup_{s,t<T,|t-s|<\delta}dist(\hat X^n_t,\hat X^n_s)$, the
following convergence holds true:
$$
\sup_n P(w_T(\hat X^n,\delta)\geq 2\eps )\to 0,\quad \delta\to 0.
$$
This allows one  to apply the above-mentioned theorem.

Our aim is to show that any limiting point of a sequence of the
distributions of the processes $\{\hat X^n\}$ gives a solution of
the martingale problem for a Walsh's Brownian motion, that was
formulated in Proposition 2. Since, by assertion 2) of Theorem 1,
this problem is correctly posed, this yields  the assertion of the
main theorem. To prove the required martingale characterization of
the limiting point, we will study the limiting behavior of
resolvents of the processes $\{\hat X^n\}$ (or, more exactly, the
Laplace transforms of their distributions; note that each of the
processes $\hat X^n$ is not Markov) in detail.

Let the vertex $\jj,$
the edge-"ray"\phantom{} with the number $\rr\leq R(\jj,\infty)$ which leaves this vertex, and the function
$\phi\in C_b([0,+\infty))$ be fixed. Consider the quantities
$$
E^{n}_{i}(t)=E(\phi(\hat X^n_t), \hat X^n_t\in
L_{\jj,\infty}^{\rr}| X^n_0=i)=E(\phi(\hat X^n_t), X^n_t\in
L_{\jj,\infty}^{\rr}| X^n_0=i),\quad t>0.
$$
By $\tau_i$, denote the  moment of the first exit of the process
$X^n$ from the neighborhood $\O_i$. For $E^{n}_{i}(\cdot)$, an
analog of the renewal equation, written at the moment $\tau_i$,
looks as
$$
E^{n}_{i}(t)=Q^{n}_{i}(t)+\sum_{k\in \P_i}\int_0^t E^{n}_{k}(t-s)
P(X^n_{\tau_i}=k,\tau_i\in ds),\eqno(6)
$$
where $Q^{n}_{i}(t)=E(\phi(\hat X^n_t), X^n_t\in
L_{\jj,\infty}^{\rr},\tau_i>t| X^n_0=i)$ (it is clear that
$Q^{n}_{i}(t)>0$ only if $i=\jj$). Considering (6) for all
$i$, we get the convolutional equation for the vector $E^{n}(t)$
composed of the components of $E_i^n(t)$. Let us introduce the Laplace
transformations
$$
U_i^n(\lambda)=\int_0^{\infty} \!\!e^{-\lambda t}
E^{n}_{i}(t)dt,\, V_i^n(\lambda)=\int_0^{\infty}\!\! e^{-\lambda
t} Q^{n}_{i}(t)dt,\, C_{i,k}^n(\lambda)= \int_0^{\infty}\!\!
e^{-\lambda t} P(X^n_{\tau_i}=k,\tau_i\in dt),
$$
and let $U^n(\lambda), V^n(\lambda),$ and $C^n(\lambda)$ be, respectively, two vectors and a
matrix composed of the components of
$U_i^n(\lambda), V_i^n(\lambda),$ and $C_{i,k}^n(\lambda)$.
Equation (6) yields
$$
U^n(\lambda)=V^n(\lambda)+C^n(\lambda)U^n(\lambda),\quad
U^n(\lambda)=[I-C^n(\lambda)]^{-1}V^n(\lambda). \eqno(7)
$$
To describe the limiting behavior of $U^n$, we need the following
lemma allowing us to  write $V^n(\lambda),C^n(\lambda)$
explicitly. Denote $\Phi(\lambda)=\int_0^\infty e^{-\lambda
t}\phi(B_t)\, dt$, where $B$ is the Brownian motion with
reflection on $[0,+\infty)$,  starting from zero.

\begin{lem}\label{Lap_calcul_lem} Let a Walsh's Brownian motion $Z$ with $M$
rays and the asymmetry parameters $p_1,\dots,p_M$ be given. Let
the points $z_1,\dots,z_m$ be marked on the rays $L_1,\dots, L_m
(m<M)$ at the distances $l_1,\dots,l_m$ from the vertex, and let
$\tau$ be the first moment when $Z$ hits one of these points.

Then, on the start of the process $Z$ from the vertex $O$,
$$
\int_0^{+\infty}\!\! e^{-\lambda t} E(\phi(Z_t), Z_t\in
L_k,\tau>t)\,dt=p_k\Phi(\lambda)\Bigl\{\sum_{j=1}^m p_j
\cth[l_j\sqrt{2\lambda}]+ \sum_{j=m+1}^M p_j\Bigr\}^{-1},\,
\lambda>0, k> m,\eqno(8)
$$
$$
\int_0^{+\infty}\!\! e^{-\lambda t} P(Z_\tau=z_k,\tau\in
dt)={p_k\over \sh[l_k\sqrt{2\lambda}]}\Bigl\{\sum_{j=1}^m p_j
\cth[l_j\sqrt{2\lambda}]+ \sum_{j=m+1}^M p_j\Bigr\}^{-1},\quad
\lambda>0, k\leq m,\eqno(9)
$$
$\cth c={\ch c\over \sh c}, \ch c={e^c+e^{-c}\over 2},\sh
c={e^c-e^{-c}\over 2}.$
\end{lem}

\emph{Proof.} It is sufficient to consider the case where, for
some $\gamma>0$, $\phi(u)=\phi(0), u\in[0,\gamma]$ (the general
case can be obtained from it by approximation of general $\phi\in
C_b([0,+\infty)$ by functions, that are constant in some
neighborhood of $0$). Let $B$ be the Brownian motion with
reflection, whose excursions have been used in the construction of
$Z$ (see Section 1.2). For $x>0$, we denote, by $\tau_x$, the time
moment of the first passage of the level $x$ by the process $B$.
Then, for $x<\min(\gamma,l_1,\dots,l_m),k>m$,
$$
E(\phi(Z_t),Z_t\in L_k,\tau>t)=\phi(0)P(Z_t\in L_k,\tau_x>t)+
$$
$$
+\sum_{j=1}^M\int_0^t E_{j,x}(\phi(Z_{t-s}),Z_{t-s}\in
L_k,\tau>t-s)P(Z_{\tau_x}\in L_j,\tau_x\in ds),
$$
where $E_{j,x}(\cdot)$ means the averaging over the distribution of the process $Z$
on the start from the point located on $L_j$ at the distance $x$ from
the vertex. We set $P(\tau_x\leq t)=T_x(t)$. By the construction of the process
$Z$, we have
$$
P(Z_t\in L_k,\tau_x>t)=p_r [1-T_x(t)],\quad  P(Z_{\tau_x}\in
L_j,\tau_x\in ds)=p_j T_x(ds).
$$
For $j\not=k$,
$$
E_{j,x}(\phi(Z_t),Z_{t}\in L_k,\tau>t)=\int_0^t
E(\phi(Z_{t-s}),Z_{t-s}\in L_k,\tau>t-s)Q_{j,x}(ds),$$ $
Q_{j,x}(s)\equiv P(W^x_{\theta_j}=0, \theta_j\leq s), $ where
$W^x$ is the Brownian motion starting from the point $x$, and
$\theta_j$ is the moment of its exit from the interval $(0,l_j)\,
(l_j\equiv +\infty$ for $j>m$). At last,
$$
E_{k,x}(\phi(Z_t),Z_{t}\in L_k,\tau>t)=\int_0^t
E(\phi(Z_t),Z_{t-s}\in L_k,\tau>t-s)Q_{j,x}(ds)+F_x(t),
$$
where $F_x(t)= E(\phi(W_t^x),\theta >t)$, $\theta$ is the moment
of exit of $W^x$ from the interval $(0,+\infty)$. Thus, we have
the convolutional equation for the function $H_k(t)=E(\phi(Z_t),
Z_t\in L_k,\tau>t)$:
$$
H_k(t)=p_k\left\{\phi(0)[1-T_x(t)]+\int_0^tF_{x}(t-s)T_x(ds)\right\}
+\sum_{j=1}^Mp_j\int_0^t\int_0^{t-s} \!\! H_r(t-s-u)
Q_{j,x}(du)T_x(ds).\eqno(10)
$$
We now introduce the Laplace transformations
$$
G_k(\lambda)=\int_0^\infty e^{-\lambda t} H_k(t)dt,\quad
S_x(\lambda)\equiv \int_0^\infty e^{-\lambda t} T_x(dt),
$$
$$
R_{j,x}(\lambda)\equiv \int_0^\infty e^{-\lambda t} Q_{j,x}(dt),
\quad \Psi_x(\lambda)=\int_0^\infty e^{-\lambda t} F_x(t)dt
$$
and rewrite  (10) in the following form:
$$
G_k(\lambda)={p_k \phi(0)\over \lambda}[1-S_x(\lambda)]+ p_k
S_x(\lambda)\Psi_x(\lambda)+\sum_{j=1}^M p_j S_x(\lambda)
R_{j,x}(\lambda) G_k(\lambda).\eqno(11)
$$
Let us differentiate (11) with respect to $x$ at the point $x=0$.
Taking into account that (see \cite{ItoMcKean}, \S 1.7)
$$
S_x(\lambda)={1\over \ch[x\sqrt{2\lambda}]},\quad
R_{j,x}(\lambda)=
\begin{cases}{\sh [(l_j-x)\sqrt{2\lambda}]\over \sh
[l_j\sqrt{2\lambda}]},& j\leq m\\
\exp[-x\sqrt{2\lambda}], & j> m \end{cases},
$$
 we have
$$
[S_x(\lambda)]'_{x=0}=0,\quad
[R_{j,x}(\lambda)]'_{x=0}=\begin{cases}-\sqrt{2\lambda} \cdot{\ch
[l_j\sqrt{2\lambda}]\over \sh
[l_j\sqrt{2\lambda}]},& j\leq m\\
-\sqrt{2\lambda}, & j> m \end{cases},
$$
and, hence,
$$
p_k [\Psi_x(\lambda)]'_{x=0}- G_k(\lambda)
\sqrt{2\lambda}\Bigl\{\sum_{j=1}^m p_j \cth[l_j\sqrt{2\lambda}]+
\sum_{j=m+1}^M p_j\Bigr\}=0. \eqno(12)
$$
Repeating literally the performed calculations for the process $B$
(which can be interpreted now as a Walsh's Brownian motion on a
graph with only one edge), we arrive at the equality
$$
[\Psi_x(\lambda)]'_{x=0}-\Phi(\lambda)
\sqrt{2\lambda}=0.\eqno(12')
$$
Relations (12) and ($12'$) yield (8). Analogously, for $H^1_k(dt)\equiv
P(Z_\tau=z_k,\tau\in dt),k\leq m$, we get
$$
H_k^1(dt)=p_k\int_0^tQ_{k,x}^1(dt-s)T_x(ds) +
\sum_{j=1}^Mp_j\int_0^t\int_0^{t-s}  H_k^1(dt-s-u)
Q_{j,x}(du)T_x(ds),
$$ where $Q_{k,x}^1(dt) \equiv
P(W^x_{\theta_k}=l_k,\theta_k\in dt)$. This equation in terms of the
Laplace transform $G_k^1(\lambda)=\int_0^\infty e^{-\lambda
t} H_r^1(dt)$ takes the form
$$
G_r^1(\lambda)=p_k S_x(\lambda)R_{k,x}^1(\lambda)+\sum_{j=1}^M p_j
S_x(\lambda) R_{j,x}(\lambda) G_r^1(\lambda),\eqno(13)
$$
where  $R_{k,x}^1(\lambda)\equiv \int_0^\infty e^{-\lambda t}
Q_{k,x}^1(dt)= {\sh [x\sqrt{2\lambda}]\over \sh
[l_k\sqrt{2\lambda}]}$. By differentiating (13) with respect to $x$ at the point $x=0$,
we get equality (9). Lemma 1 is proved.

We can suppose that $\phi(n)={1\over n}$ in (4). Since $\ch
x=1+o(x), \sh x=x+o(x^2), x\to 0,$ relation (8) yields (for a
fixed $\lambda>0$)
$$
V^n_\jj(\lambda)=\Phi(\lambda)\sqrt{2\lambda}{p^{\rr}_{\jj,\infty}\over
n}\left\{\sum_{k\in \P_\jj,r\leq R(\jj,k)} {p^{r}_{\jj,k}\over
l^{r}_{\jj,k}}\right\}^{-1}+o({1\over n}),\quad n\to +\infty,\quad
V^n_i(\lambda)=0, \quad i\not=\jj.\eqno(14)
$$
In an analogous way, relation (9) yields
$C^n(\lambda)=A^n-B^n(\lambda),$ where
$$
A_{i,k}^n={\sum\limits_{r\leq R(i,k)} {p^{r,n}_{i,k}\over
l^{r}_{i,k}}\over \sum\limits_{h\in \P_i,r\leq R(i,h)}
{p^{r,n}_{i,h}\over l^{r}_{i,h}}},\quad
B_{i,k}^n(\lambda)={\sqrt{2\lambda}\over n}{\sum\limits_{r\leq
R(i,k)} {p^{r}_{i,k}\over l^{r}_{i,k}} \sum\limits_{r\leq
R(i,\infty)} {p^{r}_{i,\infty}}\over \left\{\sum_{h\in \P_i,r\leq
R(i,h)} {p^{r}_{i,h}\over l^{r}_{i,h}}\right\}^{2}}+o({1\over
n}),\quad n\to +\infty.\eqno(15)
$$

The matrices $A^n$ converge to the matrix $A$ introduced before
the formulation of Theorem 3. This matrix is the transition
probability  matrix for a Markov chain that is ergodic by virtue
of condition (3). Hence, the matrices $A^n$ possess the same
property, beginning from some $n$. This allows us to rewrite (7)
in the form more convenient for the asymptotic analysis, by
performing a suitable change of the basis. Namely, we write the
decomposition $\Re^N=\langle e_1 \rangle +\langle \pi
\rangle^\perp, [\Re^N]^*=\langle e_1 \rangle^\perp +\langle \pi
\rangle$, where $N$ is the number of vertices in the initial graph
(knot), $\Re^N$ and $[\Re^N]^*$ are, respectively, the spaces of
column-vectors and row-vectors of dimension $N$,
$e_1=(1,\dots,1)^\intercal$, and $\pi=(\pi_1,\dots,\pi_N)$ is the
invariant distribution corresponding to $A$. Next, we choose the
basis $\{e_2,\dots,e_N\}$ in $\langle \pi \rangle^\perp$ and the
basis $\{b_2,\dots,b_N\}$ in $\langle e_1 \rangle^\perp$ in such a
way that  $\langle b_i,e_j\rangle=\delta_{ij}, i,j=1,\dots,N$,
where $\delta_{ij}$ is the Kronecker delta, $b_1\equiv \pi$, and,
for $y=(y_1,\dots,y_N)\in [\Re^N]^*$ and $x=(x_1,\dots,x_N)\in
\Re^N$, $\langle y,x\rangle\equiv \sum_{k=1}^Nx_ky_k$. Writing now
the matrix $(I-\tilde A)_{i,k}=\delta_{ik}- \langle b_i,A_{i,k}
e_k\rangle$ (i.e., we write the matrix $I-A$ w.r.t. the bases
$\{b_i\},\{e_i\}$), we get a block matrix of the form $I-\tilde A
= \begin{pmatrix} 0 & 0\\ 0 & D\end{pmatrix}$, where
$D_{i,k}=(I-\tilde A)_{i,k},i,k=2,\dots,N$. The ergodicity of the
matrix $A$ yields that
$$
[y(I-A)=0,\quad y\in \langle e_1 \rangle^\perp]\Rightarrow y=0,
\quad
[(I-A)x=0,\quad x\in \langle \pi \rangle^\perp]\Rightarrow x=0,
$$
that means that the matrix $D$ is invertible.

We perform an analogous operation for all $n$, by introducing the
bases $\{e_i^n\},\{b_i^n\}$ in such a way that $\langle
b_i^n,e_j^n\rangle=\delta_{ij},$ $e_1=(1,\dots,1)^\intercal$, and
$\pi=(\pi_1,\dots,\pi_N),$ $b_1^n=\pi^n$ being the invariant
distribution for $A^n$. It is clear that this can be performed in
such a way that $e_i^n\to e_i, b_i^n\to b_i, n\to +\infty$. Then,
writing the matrix $I-A^n$ w.r.t. the bases $\{b_i^n\},\{e_i^n\}$,
we get the matrix $I-\tilde A^n =
\begin{pmatrix} 0 & 0\\ 0
& D^n\end{pmatrix}$, where $D^n\to D,[D^n]^{-1}\to D^{-1}.$
Making now the change of variables $\tilde V_i^n(\lambda)=\langle b_i^n,V^n(\lambda)\rangle,
 \tilde U^n_i(\lambda)=\langle b_i^n,U^n(\lambda)\rangle$ in (7), we obtain
$$
\tilde U^n(\lambda)=\begin{pmatrix} \tilde B_{1,1}^n(\lambda) & \tilde B_{1,\cdot}^n(\lambda)\\
\tilde B_{\cdot,1}^n(\lambda)
& D^n+\tilde B_{\cdot,\cdot}^n(\lambda)\end{pmatrix}^{-1} \tilde V^n(\lambda),\eqno(16)
$$
where $\tilde B_{i,k}^n(\lambda)=\langle b_i^n,B^n(\lambda)
e_k^n\rangle$, $\tilde B_{\cdot,\cdot}^n(\lambda)=(\tilde
B_{i,k}^n(\lambda))_{i,k=2}^N$, $\tilde
B_{1,\cdot}^n(\lambda)=(\tilde B_{1,2}^n(\lambda),\dots,\tilde
B_{1,N}(\lambda)),$ $\tilde B_{\cdot,1}^n(\lambda)=(\tilde
B_{2,1}^n(\lambda),\dots,\tilde B_{N,1}(\lambda))^\intercal.$ To
invert the block matrix in (16), we will use a suitable
modification of the Gauss method. We note that, at sufficiently
great $n$, the matrix $D^n+\tilde B_{\cdot,\cdot}^n(\lambda)$ is
invertible. Next, by virtue of (15), $\tilde
B_{i,k}^n(\lambda)=O({1\over n})$ and
$$
n\tilde B_{1,1}^n(\lambda)\to
\sqrt{2\lambda}\sum_{i=1}^N\left\{{\pi_i \sum\limits_{r\leq
R(i,\infty)} {p^{r}_{i,\infty}}\over \sum_{h\in \P_i,r\leq R(i,h)}
{p^{r}_{i,h}\over l^{r}_{i,h}}}\right\}>0,\quad  n\to+\infty.
$$
Therefore, the number $\gamma_n(\lambda)= \tilde
B_{1,1}^n(\lambda)- \tilde B_{1,\cdot}^n(\lambda) [D^n+\tilde
B_{\cdot,\cdot}^n(\lambda)]^{-1}\tilde B_{\cdot,1}^n(\lambda)$ is
nonzero for great $n$. Then,  denoting $\tilde
V_{\cdot}^n(\lambda)=(\tilde V_{2}^n(\lambda),\dots,\tilde
V_{N}(\lambda))^\intercal,$ $\tilde U_{\cdot}^n(\lambda)=(\tilde
U_{2}^n(\lambda),\dots,\tilde U_{N}(\lambda))^\intercal,$ we get
(see \cite{Heinrich}, Lemma 2.3, and, in particular, formula (12))
that, for such $n$,
$$
\tilde U_{1}^n(\lambda)= [\gamma_n(\lambda)]^{-1}\Bigl(\tilde V_1^n(\lambda)-
\langle \tilde B_{1,\cdot}^n(\lambda),[D^n+\tilde B_{\cdot,\cdot}^n(\lambda)]^{-1}\tilde V_{\cdot}^n(\lambda)\rangle
\Bigr),
$$
$$
\tilde U_{\cdot}^n(\lambda)=[D^n+\tilde B_{\cdot,\cdot}^n(\lambda)]^{-1}
(-\tilde U_{1}^n(\lambda)\tilde B_{\cdot,1}^n(\lambda)+\tilde V_{\cdot}^n(\lambda)).
$$
These relations and the relations
$$
\tilde V_{i}^n(\lambda)=O({1\over n}),\quad n\tilde
V_{1}^n(\lambda)\to \Phi(\lambda) \sqrt{2\lambda}\cdot \pi_\jj
p^{\rr}_{\jj,\infty} \left\{\sum_{k\in \P_\jj,r\leq R(\jj,k)}
{p^{r}_{\jj,k}\over l^{r}_{\jj,k}}\right\}^{-1},\quad n\to
+\infty,
$$
that are valid by virtue of (14), imply that, as $n\to +\infty$,
$$
\tilde U_{1}^n(\lambda)\to \Phi(\lambda) {{\pi_\jj\over \beta_\jj}
p^{\rr}_{\jj,\infty} \over \sum_{i=1}^N \sum_{r\leq R(i,\infty)}
{\pi_i\over \beta_i} p^{r}_{i,\infty}}, \quad \tilde
U_{i}^n(\lambda)\to 0, \quad i=2,\dots,N.
$$
We recall that $\beta_i=\sum_{k\in \P_i,r\leq R(i,k)}
{p^{r}_{i,k}\over l^{r}_{i,k}}$. Performing the inverse change of
the variables $\tilde U\to U$, we get that, for an arbitrary
vertex $i$,
$$
U_{i}^n(\lambda)=\int_0^{\infty} \!\!e^{-\lambda t} E(\phi(\hat
X^n_t), \hat X^n_t\in L_{\jj,\infty}^{\rr}| X^n_0=i) dt \to
\Phi(\lambda)\cdot \hat p^\rr_{\jj,\infty}, \eqno(17)
$$
where $\hat p^\rr_{\jj,\infty}$ is given by equality (5). For an
arbitrary point $x\in \G$,
$$
E(\phi(\hat X^n_t), \hat X^n_t\in L_{\jj,\infty}^{\rr}| X^n_0=x)=
E(\phi(\hat X^n_t), \tau>t|  X^n_0=x)\1_{x\in
L_{\jj,\infty}^{\rr}}+
$$
$$
+\sum_{i=1}^N\int_0^t E(\phi(\hat X^n_{t-s}), \hat X^n_{t-s}\in
L_{\jj,\infty}^{\rr}| X^n_0=i)P(\tau\in ds, X^n_\tau=i), \eqno(18)
$$
where $\tau$ is the moment when the process $X^n$ hits  one of the
vertices. For $x$ lying on one of the rays, the distribution of
$\tau$ coincides with the distribution of the same moment for a
Walsh's Brownian motion (or a Wiener process on a half-line). On
the other hand, it is easy to verify that $\tau\toP 0$ for $x\in
\G_0$ ($\G_0=\G\backslash \hat \G$ is the joint of all the edges
contracting into a point), and this convergence is uniform in
$x\in \G_0$. Then, for an arbitrary continuous function $\psi$ on
the limiting bundle of rays $\hat \G$, relations (17) and (18)
yield
$$
\int_0^{\infty} \!\!e^{-\lambda t} E(\psi(\hat X^n_t))| X^n_0=x)
dt \to R_\lambda\psi(O),\quad x\in \G_0,\eqno(19)
$$
$$
\int_0^{\infty} \!\!e^{-\lambda t} E(\psi(\hat X^n_t))| X^n_0=x)
dt \to R_\lambda\psi(x), \quad x\in \hat \G,\eqno(19')
$$
where $R_\lambda\psi(x)=\int_0^{\infty} \!\!e^{-\lambda t}
E(\psi(\hat X_t))| X_0=x)$ is the resolvent of the Walsh's
Brownian motion $\hat X$, whose asymmetry parameters are given by
equality (5). Moreover, the convergence in (19) and ($19'$) is
uniform on $\G_0$ and every bounded subset $\hat\G$, respectively.

We now can present, at last, a martingale characterization of the
limiting points of the sequence $\{\hat X^n\}$. Let $\tilde
 X$ be one of such points which is the limit of the subsequence
$\{\hat X^{n_k}\}$. By $A$, we denote the infinitesimal operator
of the above-mentioned Walsh's Brownian motion $\hat X$, whose
domain has been described in Proposition 2. For arbitrary
$\lambda>0, \phi\in {\Deu}_A$, $t_0>0,{t_1,\dots,t_m}\subset
[0,t_0], G\in C_b(\Re^m)$, by using the Markov property of the
processes $\hat X^{n_k}$ and relations (19) and ($19'$) for
$\psi=\phi$ and $\psi=A\phi$, we get the equality
$$
E\int_{t_0}^\infty \!\!e^{-\lambda (t-t_0)}  G(\tilde
X_{t_1},\dots,\tilde X_{t_m})[\lambda \phi(\tilde
X_t)-A\phi(\tilde X_t)]dt=
$$
$$
=\lim_{k\to +\infty}E\int_{t_0}^\infty \!\!e^{-\lambda (t-t_0)}
G(\hat X_{t_1}^{n_k},\dots,\hat X_{t_m}^{n_k})[\lambda \phi(\hat
X_{t}^{n_k})-A\phi(\hat X_{t}^{n_k})]dt
$$
$$
=E  G(\tilde X_{t_1},\dots,\tilde X_{t_m})[\lambda
R_\lambda\phi(\tilde X_{t_0}) -R_\lambda A\phi(\tilde X_{t_0})]=E
G(\tilde X_{t_1},\dots,\tilde X_{t_m})\phi(\tilde X_{t_0}).
\eqno(20)
$$
We now repeat the arguments given in \cite{FreWen}, Section 2: relation (20)
yields
$$
\int_{t_0}^\infty \!\!\lambda e^{-\lambda t} E G(\tilde
X_{t_1},\dots,\tilde X_{t_m})\left [\phi(\tilde X_t)- \phi(\tilde
X_{t_0})- \int_{t_0}^t A\phi(\tilde X_s)\, ds\right]\,dt=0.
$$
Since a continuous function is uniquely determined by its Laplace
transform, we have that $ E G(\tilde X_{t_1},\dots,\tilde X_{t_m})
[\phi(\tilde X_t)- \phi(\tilde X_{t_0})- \int_{t_0}^t A\phi(\tilde
X_s)\, ds]=0.$ Since $G$ and $t_1,\dots t_m$ are arbitrary, this
implies that the process $\phi(\tilde X_t)- \int_{0}^t
A\phi(\tilde X_s)\, ds$ is a martingale. That is, $\tilde X$ is a
solution of the martingale problem presented in Proposition 2. In
addition, the distribution of $\tilde X_0$ is known and coincides
with $\mu$. Then, by virtue of Theorem 1, the distribution of
$\tilde X$ coincides with the distribution of a Walsh's Brownian
motion with the initial distribution $\mu$. Because we took, as
$\tilde X$, an arbitrary limiting point of the sequence of
distributions $\{\hat X^n\}$, this means that the whole sequence
converges to the indicated distribution. Theorem 3 is proved.

\vskip 20pt \centerline{\textsc{3. Examples}} \vskip 20pt

In this section, we present two examples of the application of
Theorem 3. First, consider the example, in which the limiting
graph $\hat \G$ consists of one vertex and two half-lines. In this
case, the limiting process is the well known skew Brownian motion
(see \cite{ItoMcKean}, \cite{HarShepp}).

{\bf Example 1.} Let $\G$ consist of two vertices $\{1,2\}$ and
four edges $L_{1,2}, L_{2,1}, L_{1,\infty},$ and $L_{2,\infty}$.
Let the lengths of the edges $L_{1,2}$ and $L_{2,1}$ be $
l_{1,2}^n=n^{-1}\theta$ and $l_{2,1}^n=n^{-1}\kappa,$ where
$\theta,\kappa >0$ are some constants. We suppose that the drift
coefficient is zero and the  diffusion coefficient  is unity on
all the edges. We write the asymmetry parameters for the vertices
$\{1,2\}$ in the form
$$
p_{1,\infty}^n={1+q_1^n\over 2},\quad p_{1,2}^n={1-q_1^n\over
2},\quad p_{2,\infty}^n={1+q_2^n\over 2},\quad
p_{2,1}^n={1-q_1^n\over 2},\quad q_i^n\to q_i\in(-1,1),\eqno(21)
$$
$i=1,2.$ Then, by virtue of Theorem 3 (the calculations are simple
and omitted), the limiting process possesses the asymmetry
parameters $p_{\pm}={1\pm q\over 2}$, where
$$
{1+q\over 1-q}= {\theta\over \kappa}\cdot  {1+q_1\over 1-q_1}\cdot
{1-q_2\over 1+q_2}.\eqno(22)
$$
Consider two special cases. First, if $l_{1,2}^n=l_{2,1}^n$, we
have, in essence, the graph with one non-oriented edge contracting
into a point or, from the other viewpoint, the diffusion process
with two singular points ("membranes") contracting into a single
point. We can rewrite now formula (22) as
$q=\th(\lambda_1-\lambda_2),$ $\lambda_{i}=\arcth q_i$, which
reproduces the result mentioned in Introduction. The difference in
signs is caused by the fact that the notation of the asymmetry
parameters used in \cite{Zaitseva} differs from (21).

Another interesting special case arises if $q_i^n\equiv 0$. In
this situation, the prelimiting process is a mixture of two
Brownian motions on the half-lines $(-\infty, n^{-1} \theta],
[-n^{-1} \kappa, +\infty)$, for which the switch from one motion
to another one happens at the moments of hitting of the point
$n^{-1} \theta$ (for the first motion) or $-n^{-1} \kappa$ (for
the second one). After the switch, a new motion starts from the
point $0$. Such a process is naturally interpreted as \emph{the
Brownian motion with buffer zones} on a straight line.

It is known that the asymmetry of a skew Brown motion at zero can
be interpreted in different ways: as the result of "a gambling of
excursions"\phantom{} with unequal probabilities
(\cite{ItoMcKean}), the presence of a singular drift coefficient
$a=q\delta_0$ (\cite{HarShepp}), or, what is close to the
previous,  the presence of a semipermeable membrane at the point 0
(\cite{portenko}). Theorem 3 allows us to propose one more,
apparently completely new interpretation of the skew Brownian
motion. Let us consider the ordinary Brownian motion (without any
singularity!) and introduce the buffer zones in the neighborhood
of zero with a fixed ratio of the lengths of the zones
$\vartheta={\theta\over \kappa}$. Then, "from the macroscopic
viewpoint"\phantom{} (i.e., when the size of the zones tends to
zero), this motion has the form of a skew Brownian motion with the
parameter $q={\vartheta-1\over \vartheta+1}={\theta-\kappa\over
\theta+\kappa}$ (see Fig. 1).

\begin{center}
\includegraphics{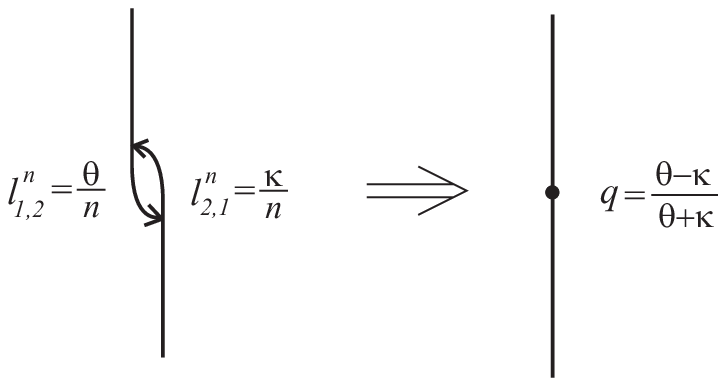}%

Fig. 1.
\end{center}

{\bf Example 2.}  Let $\hat \G$ be a bundle consisting of $m$
half-lines, and let the diffusion process $\hat X$ on $\hat \G$
have the coefficients of drift and diffusion $a^{r}$ and
$\sigma^r, r=1,\dots, m$ and the asymmetry parameters $p^r,
r=1,\dots, m$. Additionally, we assume that the functions $a^{r},
\sigma^r, r=1,\dots, m$ satisfy the local Lipschitz condition on
$\ax$, so that, from the start from some point on the edge till
the moment of hitting of the vertex, the process $\hat X$ can be
represented as the strong solution of an SDE driven by a
one-dimensional Wiener process. At the same time, if at least
three numbers from the collection  $\{p^r\}$ are nonzero, then the
vertex of the graph $\hat \G$ is a triple point, and the diffusion
$\hat X$ contains the Tsirelson's singularity and, in particular,
cannot be presented  as the strong solution to any system of SDEs
driven by a (multidimensional or even infinitedimensional) Wiener
process.

We now construct an approximation of this diffusion (in the sense
of the convergence in distribution) by the sequence of diffusions
not containing such a singularity. Consider the graph $\G$
consisting of $m$ vertices enumerated by the numbers $1,\dots, m$
and $2m$ edges. One edge-"ray" and one edge-"segment" leave each
vertex. The edge-"segment"\phantom{} joins the $i$-th vertex with
the $i+1$-th one (to shorten the notation, we use the agreement
$m+1\equiv 1$). We put
$$
a_{i, \infty}^n=a^i,\quad  \s_{i,\infty}^n=\s^i,\quad
a_{i,i+1}^m\equiv a^i(0),\quad  \s_{i,i+1}^n\equiv
\s^i(0),\eqno(23)
$$
$$
p_{i,\infty}^n= p_{i,i+1}^n={1\over 2},\quad l_{i,i+1}^n={1\over
n}p^i, \eqno(24)
$$
$i=1,\dots, m$ (we drop the superscript $r$, because each two
vertices are connected by at most one edge, and exactly one
edge-"ray" leaves each vertex). Let us consider a sequence of
diffusion processes $\{X^n\}$ on $\G$ with parameters (23),(24)
and apply Theorem 3 to it. The limiting graph has the single
vertex, therefore, we drop the superscript $\hat i$ in the
notations given below. We have $\alpha_{i,i+1}=1$ and
$\alpha_{i,j}=0$ for $j\not=i+1$. Therefore,
$$
 \beta_i=1,\quad A_{i,j}=\begin{cases}1,& j=i+1,\\ 0& j\not=i+1\end{cases},\quad
\pi_i={1\over m}, \quad i={1\dots, m}.
$$
Then, as far as $\sum_{r} p^r=1$, the normalizing constant $P$ in
formula (5) is $P={1\over m}$. Hence,
the limiting asymmetry parameters are equal to $p^1,\dots, p^m$. The coefficients
$a_{i,j}^n, \s_{i,j}^n$, satisfy, obviously, the conditions of Theorem
2. Hence, the projections $\hat X^n$ of the processes $X^n$ on $\hat \G$
converge in distribution in $C(\ax, \hat \G)$ to the process $\hat X$
(see Fig. 2).

\begin{center}
\includegraphics{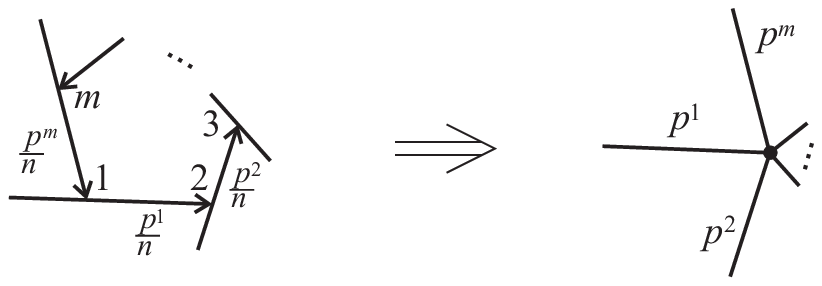}%

Fig. 2.
\end{center}

By construction, the process $X^n$ has no triple points. Let us
show that this process can be presented as a mixture (in the same
sense as that in Example 1) of the collection of strong solutions
of one-dimensional SDEs. This means that $X^n$ can be obtained
from a  one-dimensional Wiener process  "in real time".

Consider $m$ diffusion processes $Y_1^n, \dots, Y^n_m$, and let
the $i$-th process be defined on a half-line $L_i^n=[-{1\over
n}p^i, +\infty)$ and have the coefficients
$$
a_i(x)=a^i(0)\1_{[-{1\over n}p^i, 0]}(x)+a^i(x)\1_\ax(x), \quad
\s_i(x)=\s^i(0)\1_{[-{1\over n}p^i, 0]}(x)+\s^i(x)\1_\ax(x).
$$
By virtue of the Lipschitz property of the coefficients, each of
the processes $Y_i^n$ can be represented as the strong solution of
a one-dimensional SDE prior to its exit on the boundary of a
half-line. The process $X^n$ can be represented now as follows: at
every time moment, the point $X^n$ is located on one of the
half-lines $L_i^n$ and moves along it by the law defined by the
process $Y_i^n$ till the moment of its exit on the boundary of a
half-line. Then the point passes to the point with the coordinate
$0$,  located on the half-line $L_{i+1}^n$, and proceed moving
along $L_{i+1}^n$ via the law of $Y_{i+1}^n$, etc. The presented
interpretation implies that the filtration generated by the
process $X^n$ can be obtained by  a morphism from the filtration
generated by a one-dimensional Wiener process. That is, the
process $X^n$ does not contain Tsirelson's singularity.

 Let us summarize: "untwisting"\phantom{} a vertex
with multiplicity $\geq 3$, i.e., representing a vertex with
multiplicity $\geq 3$ as a result of the "contraction into a
point"\phantom{} of a knot, all the vertices of which have
multiplicity $\leq 2$, we have obtained the approximation of a
process containing the Tsirelson's singularity by processes
without such a singularity. It is clear that one can apply the
same trick to graphs of arbitrary configuration,
"untwisting"\phantom{} each vertex with multiplicity $\geq 3$.

\end{document}